\documentclass{llncs}

\usepackage{url}
\usepackage{cite}
\usepackage{graphicx}
\usepackage{framed}

\begin{document}

\title{No, This is not a Circle!}
\author{Zolt\'an Kov\'acs}
\institute{
The Private University College of Education of the Diocese of Linz\\
Salesianumweg 3, A-4020 Linz, Austria\\
\email{zoltan@geogebra.org}
}

\maketitle              

\begin{center}
\includegraphics[scale=0.4]{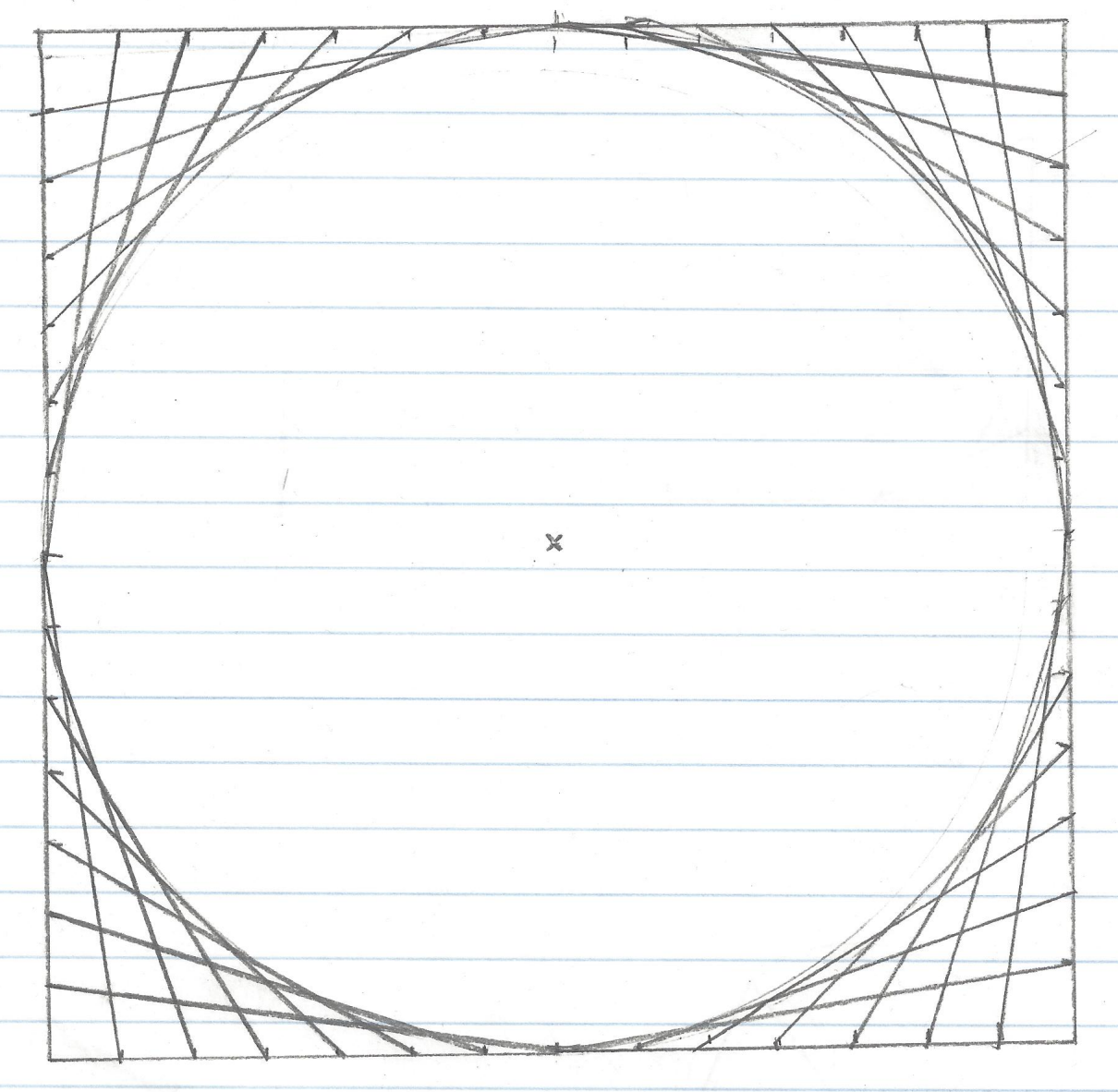}
\end{center}

\begin{abstract}

A popular curve shown in introductory maths textbooks, seems like a circle.
But it is actually a different curve. This paper discusses some elementary
approaches to identify the geometric object, including novel technological means by using GeoGebra.
We demonstrate two ways to refute the false impression, two suggestions to find a correct
conjecture, and four ways to confirm the result by proving it rigorously.

All of the discussed approaches can be introduced in classrooms at various levels
from middle school to high school.

\keywords{string art, envelope, GeoGebra, computer algebra, computer aided mathematics education, automated theorem proving}

\end{abstract}
\section{But it looks like a circle}

One possible anti boredom activity is to simulate string art in a chequered notebook,
as seen below the title of this paper.
This kind of activity is easy enough to do it very early, even as a child during
the early school years. The resulting curve, the contour of the ``strings'', or more
precisely, a curve whose tangents are the strings, is called an \textit{envelope}.

\begin{figure}
\begin{center}
\includegraphics[width=0.4\textwidth]{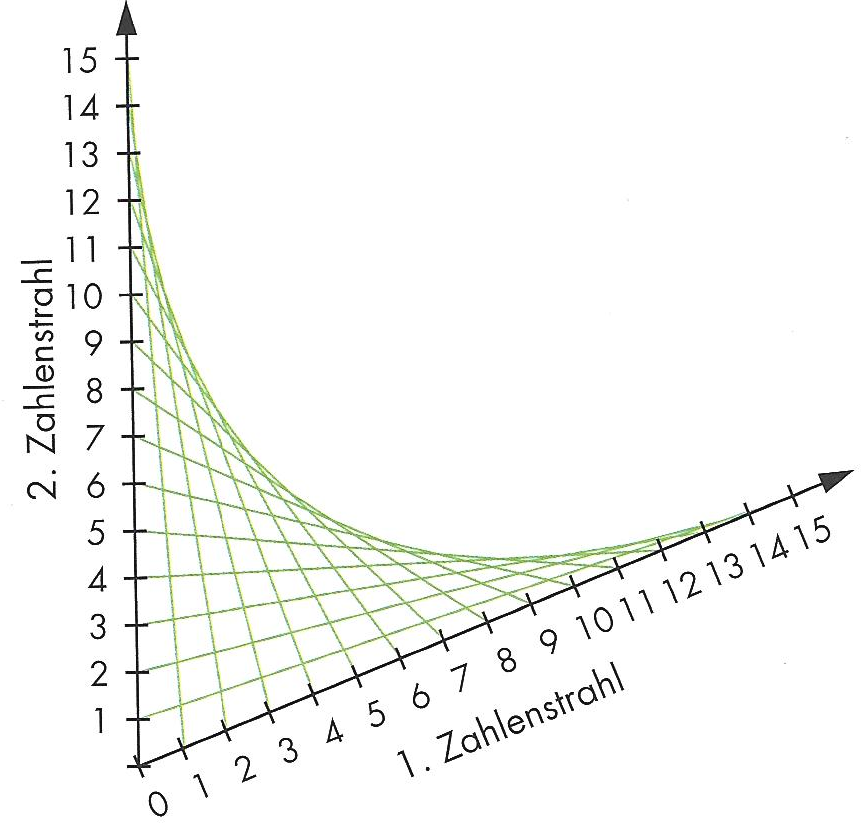}
\includegraphics[width=0.5\textwidth]{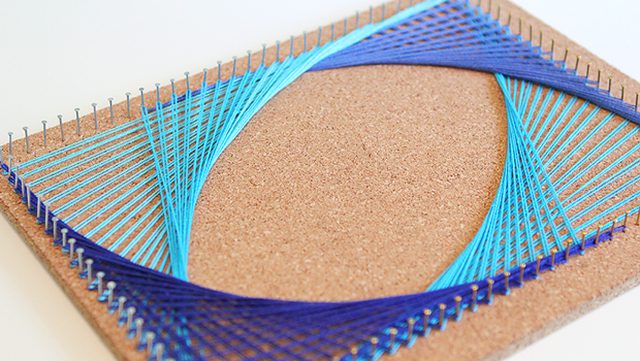}
\caption{An activity for young learners on the left: Join the numbers on each ray by
a segment to make sums 16 \cite{mathematix}. Such activities are also called `string art' when
they are performed by sewing a thread on some fabric or other material. On the right
there is a combinination of 4 simple constructions, but produced in a different layout than
the one shown at the beginning of this paper \cite{Morgan2013}.}
\label{fig:mathematix}
\end{center}
\end{figure}

According to Wikipedia \cite{wiki:envelope}, \textit{an envelope of a family of curves in the plane} is
a curve that is tangent to each member of the family at some point.\footnote{This definition is however polysemic:
the Wikipedia page lists other non-equivalent ways to introduce the notion of envelopes. See \cite{BR-survey-envelopes}
for a more detailed analysis on the various definitions.} Let
us assume that the investigated envelope---below the title of this paper---which is defined similarly as
the learner activity in Fig.~\ref{fig:mathematix}, is a circle. In the investigated
envelope it will be assumed that a combination of 4 simple constructions is used, the axes are perpendicular, and the sums of the joined numbers
are 8. To be more general, these sums may be changed to different (but fixed) numbers.
These sums will be denoted by $d$ to recall the \textit{distance} of the origin and the furthermost point
for the exterior strings.

\begin{figure}
\begin{center}
\includegraphics[scale=0.25]{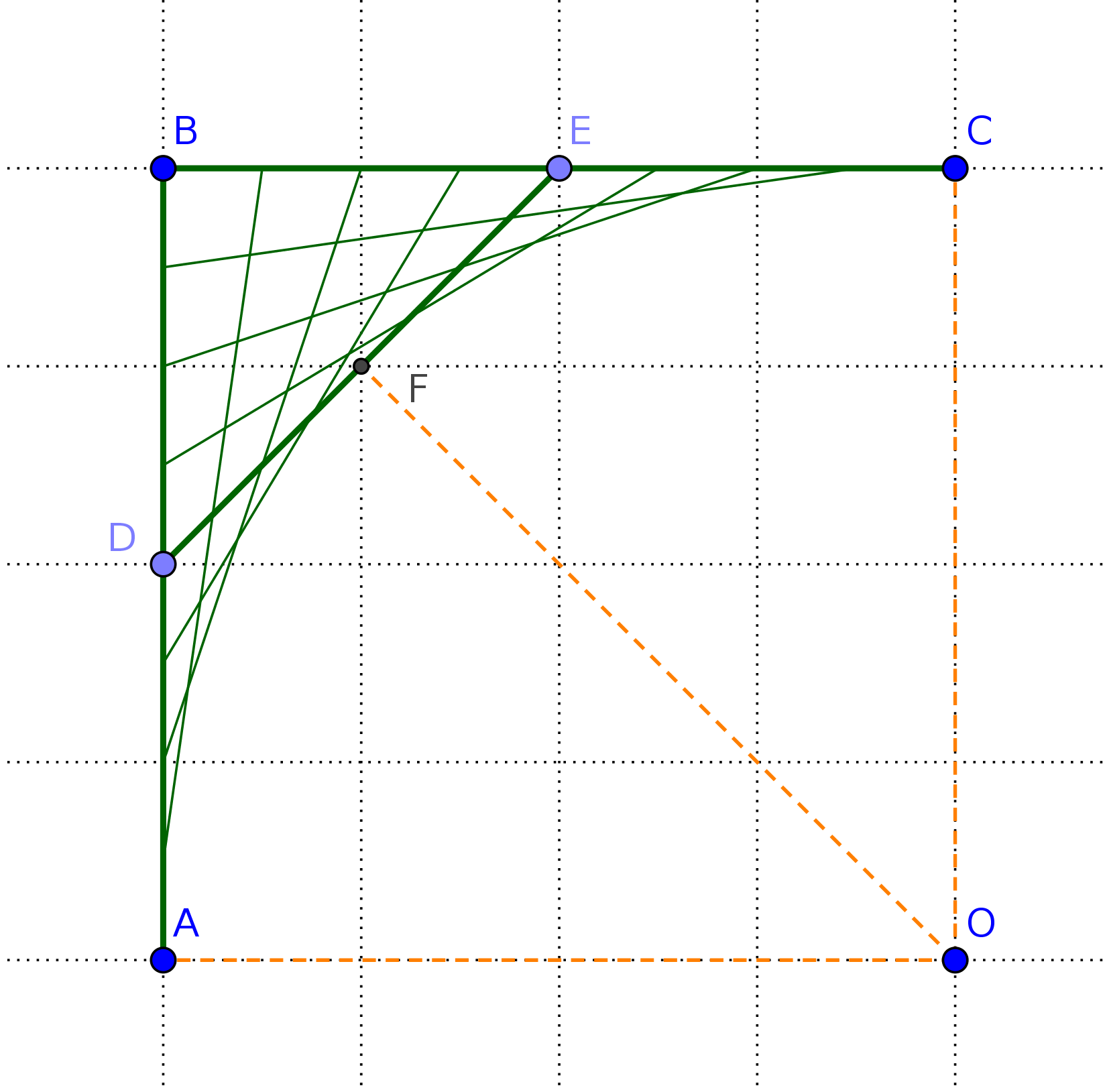}
\caption{Considering three strings from the family and their distance from the assumed center}
\label{fig:topleft}
\end{center}
\end{figure}

By using the \textit{assumption} of the circle property, in our case the
family of the strings must be equally far from the center of the circle. This needs to
be true because the circle is the only curve whose tangents are equally far from the center. Due
to symmetry of the 4 parts of the figure, the only possible center for the circle is the
midpoint of the figure. Let us consider the top-left part of the investigated figure
(Fig.~\ref{fig:topleft}). On the left and the top the strings $AB$ and $BC$ have the
distance $d=OA=OC$ from center $O$. On the other hand, the diagonal string $DE$
has distance $OF=\frac{3}{4}\cdot d\cdot\sqrt2$ from the assumed center,
according to the Pythagorean theorem.
This latter distance is approximately $1.06\cdot d$, that is, more than $d$.
Consequently, the curve cannot be an exact circle. That is, it is indeed
not a circle.

In schools the Pythagorean theorem is usually introduced much later than the
students are ready to simply measure the length of $OA$ and $OF$ by using
a ruler. The students need to draw, however, a large enough figure because
the difference between $OA$ and $OF$ is just about 6\%. Actually, both
methods obviously prove that the curve is different from a circle, and the
latter one can already be discussed at the beginning of the middle school.

\section{OK, it is not a circle---but what is it then?}

\begin{figure}
\begin{center}
\includegraphics[scale=0.7]{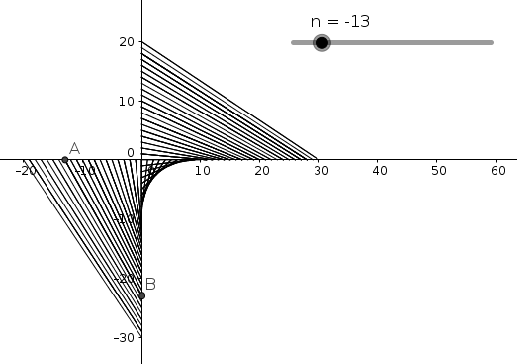}
\caption{Let us assume that $d=10$ and create a GeoGebra applet as seen in the figure.
(Actually, an arbitrary $d$ can be chosen without loss of generality.)
Now slider $n$ in range $[-20,30]$ with integer values creates points $A=(n,0)$ and $B=(0,n-10)$.
The family of segments $AB$ may enlighten which curve is the investigated one.}
\label{fig:extended}
\end{center}
\end{figure}

Let us continue with a possible classroom solution of the problem. Since
the strings are easier to observe than the envelope, it seems logical
to collect more information about the \textit{strings}. Extending the definition of
the investigated envelope by continuing the strings to both directions,
we learn how the slope of the strings changes while continuing the
extension more and more (Fig.~\ref{fig:extended}). Here we remark that the top-left
part of the investigated envelope is now mirrored about the first axis (cf.~Fig.~\ref{fig:mathematix}),
therefore not the sums but the \textit{differences} will be constant, namely $d=10$.

The strings in the extension support the idea that the tangents of the curve, when
$|n|$ is large enough, are almost parallel to the line $y=-x$. This observation
may refute the opinion that the curve is eventually a hyperbola (which has two
asymptotes, but they are never parallel).

\begin{figure}
\begin{center}
\includegraphics[scale=0.7]{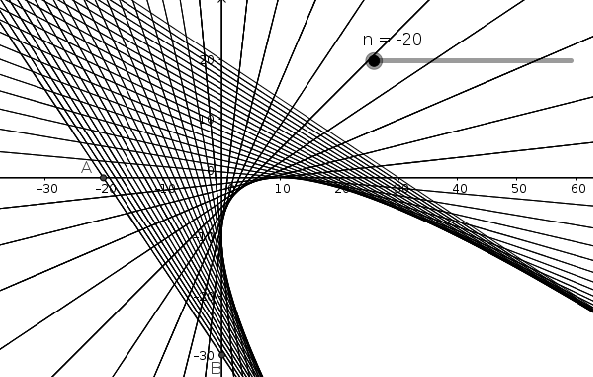}
\caption{Instead of segments as in Fig.~\ref{fig:extended} we use lines}
\label{fig:extended2}
\end{center}
\end{figure}

On the other hand, by changing the segments in Fig.~\ref{fig:extended} to lines
an obvious conjecture can be claimed, that is, the curve is a parabola (Fig.~\ref{fig:extended2}).
Thus the observed curve must be a union of 4 parabolic arcs.

\section{We have a conjecture---can we verify that?}

\begin{figure}
\begin{center}
\includegraphics[scale=0.35]{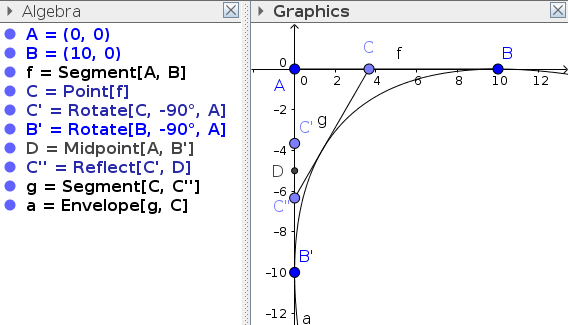}
\caption{A GeoGebra applet to compute and plot the parabola. Thanks to Michel Iroir and
No\"el Lambert for suggesting this construction method. A similar approach can be
found at \protect\url{http://dev.geogebra.org/trac/browser/trunk/geogebra/test/scripts/benchmark/art-plotter/tests/string-art-simple.ggb}.
}
\label{fig:parabola}
\end{center}
\end{figure}

A GeoGebra applet in Fig.~\ref{fig:parabola} can explicitly compute the equation of the
envelope and plot it accurately. (See \cite{BR-survey-envelopes} for a detailed survey
on the currently available software tools to visualize envelopes dynamically.) For technical reasons a slider cannot be used
in this case---instead a purely Euclidean construction is required as shown in the figure.
Free points $A$ and $B$ are defined to set the initial parameters of the applet,
and finally segment $g=CC''$ describes the family of strings. The command \texttt{Envelope[$g$,$C$]}
will then produce an implicit curve, which is in this concrete case $x^2+2xy-20x+y^2+20y=-100$.

GeoGebra uses heavy symbolic computations in the background to find this curve \cite{mcs}. Since
they are effectively done, the user may even drag points $A$ and $B$ to different positions
and investigate the equation of the implicit curve. They are recomputed quickly enough to
have an overview on the resulted curve in general---they are clearly quadratic algebraic
curves in variables $x$ and $y$.

Without any deeper knowledge of the classification of algebraic curves, of course, young
learners cannot really decide whether the resulted curve is indeed a parabola. Advanced learners
and maths teachers could however know that all real quadratic curves are either circles, ellipses,
hyperbolas, parabolas, a union of two lines or a point in the plane. As in the above, we can
argue that the position of the strings as tangents support only the case of parabolas here.

\begin{figure}
\begin{center}
\includegraphics[scale=0.25]{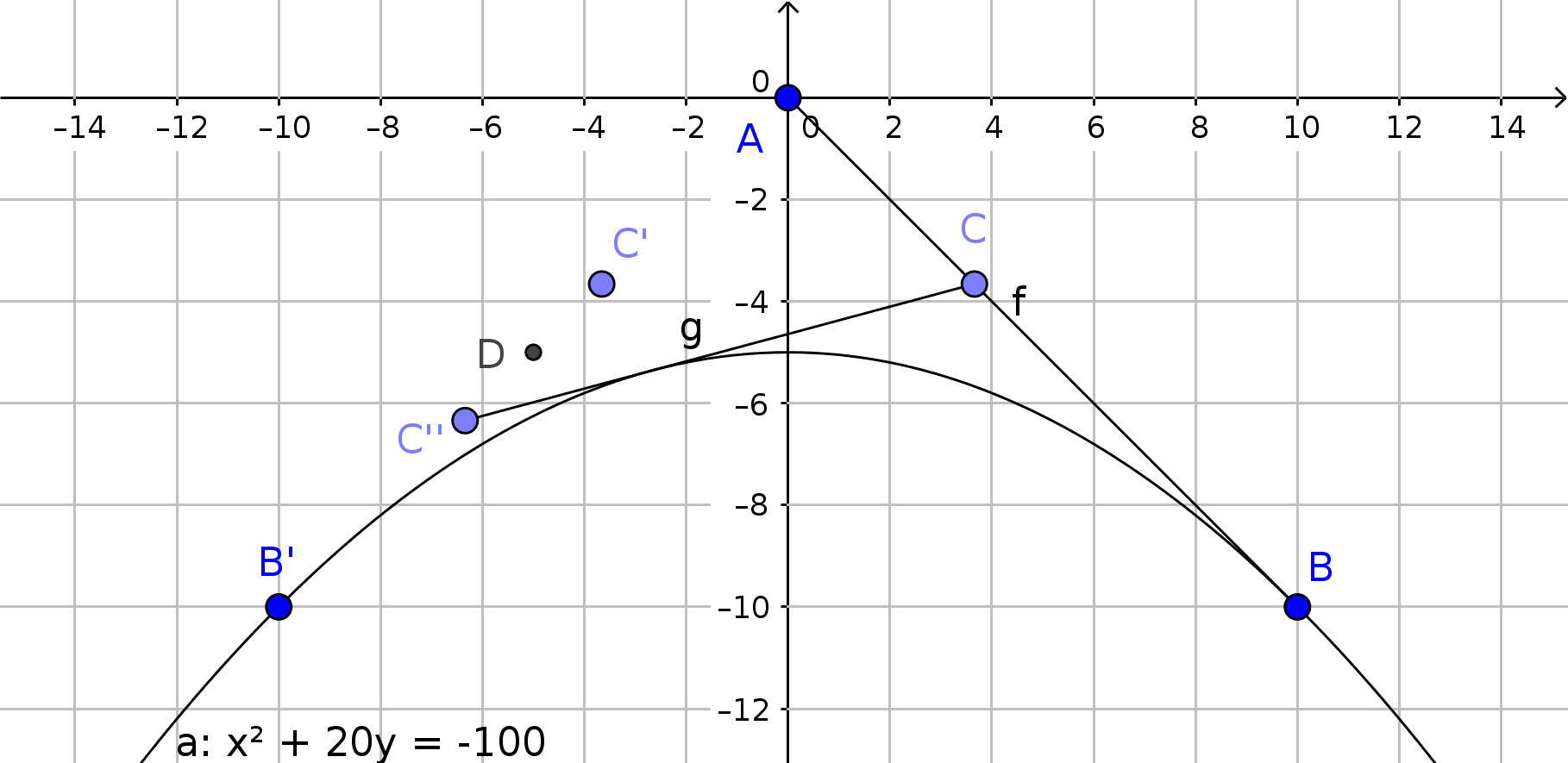}
\caption{By choosing $B=(10,-10)$ we obtain a simpler equation for the implicit curve}
\label{fig:parabola2}
\end{center}
\end{figure}

On the other hand, for young learners we can still find better positions for $A$ and $B$.
It seems quite obvious that the curve remains definitely the same (up to similarity), so
it is a free choice to define the positions of $A$ and $B$. By keeping $A$ in the origin
and putting $B$ on the line $y=-x$ we can observe that the parabola is in the form
$y=ax^2+bx+c$ which is the usual way how a parabola is introduced in the classroom. 
(In our case actually $b=0$.) For example,
when $B=(10,-10)$, the implicit curve is $x^2+20y=-100$, and this can be easily converted
to $y=-\frac{1}{20}x^2-5$ (Fig.~\ref{fig:parabola2}).

This result is computed by using precise algebraic steps in GeoGebra. One can check these
steps by examining the internal log---in this case 16 variables and 11 equations will be
used including computing a Jacobi determinant and a Gr\"obner basis when eliminating
all but two variables from the equation system. That is, GeoGebra actually provides
a proof, albeit its steps remain hidden for the user. (Showing the detailed steps when
manipulating on an equation system with so many variables makes no real sense from
the educational point of view: the steps are rather mechanical and may fill hundreds of pages.)

As a conclusion, it is actually proven that the curve is a parabola. Of course, learners
may want to \textit{understand why} it that curve is.

\section{Proof in the classroom}

Here we provide two simple proofs on the fact
that the envelope is a parabola in Fig.~\ref{fig:parabola2}. 
The first method follows \cite{LociEnvelopes-arXiv}.

We need to
prove that segment $CC''$ is always a tangent
of the function $y=-\frac{1}{20}x^2-5$. First we compute the
equation of line $CC''$ to find the intersection
point $T$ of $CC''$ and the parabola.

We recognize that if point $C=(e,-e)$,
then point $C''=(e-10,e-10)$.

Now we have two possible approaches to continue.

\begin{enumerate}
\item
Since line $CC''$ has an
equation in form $y=ax+b$, we can set up equations for points
$C$ and $C''$ as follows:
\begin{equation} \label{eq1}
-e = a\cdot e+b
\end{equation}
and
\begin{equation} \label{eq2}
e-10 = a\cdot(e-10)+b.
\end{equation}
Now (\ref{eq1})$-$(\ref{eq2}) results in $a=1-\frac{1}{5}e$ and thus, by using
$(1)$ again we get $b=-2e+\frac{1}{5}e^2$.

Second, to obtain intersection point $T$ we consider equation
$ax+b=-\frac{1}{20}x^2-5$ which can be reformulated to search the
roots of quadratic function $\frac{1}{20}x^2+ax+b+5$. If and only
if the discriminant of this quadratic expression is zero, then
$CC''$ is a tangent. Indeed, the determinant is
$a^2-4\cdot\frac{1}{20}\cdot(b+5)=a^2-\frac{b}{5}-1$ which is, after
inserting $a$ and $b$, obviously zero.

\item
Another method to show that $CC''$ is a tangent of the parabola is to
use elementary calculus. School curricula usually includes computing
tangents of polynomials of the second degree.

Let $T=(t,-\frac{1}{20}t^2-5)$. Now the steepness of the tangent of the parabola
in $T$ is $(-\frac{1}{20}t^2-5)'=-\frac{1}{10}t$. It means that the equation
of the tangent is $y=-\frac{1}{10}tx+b$, here $b$ can be computed by using
$x=t$ and $y=-\frac{1}{20}t^2-5$, that is $b=\frac{1}{20}t^2-5$. The equation
of the tangent is consequently
\begin{equation} \label{eq3}
y = -\frac{1}{10}tx+\frac{1}{20}t^2-5.
\end{equation}

Let us assume now that $C$ and $C''$ are the intersections of the tangent
and the lines $y=-x$ and $y=x$, respectively. The $x$-coordinate of $C$
can be found by putting $y=-x$ in (\ref{eq3}), it is
$$x_C=\frac{\frac{1}{20}t^2-5}{\frac{1}{10}t-1}.$$ On the other hand, the
$x$-coordinate of $C''$ can be found by putting $y=x$ in (\ref{eq3}), it
is $$x_{C''}=\frac{\frac{1}{20}t^2-5}{\frac{1}{10}t+1}.$$ By using some
basic algebra it can be confirmed that $x_C-10=x_{C''}$, that is $CC''$
is indeed a string.

\end{enumerate}

The second proof is technically longer than the first one but still achievable
in many classrooms.

Both approaches are purely analytical proofs without any knowledge
of the synthetic definition of a parabola. The fact is that in many classrooms,
unfortunately, the synthetic definition is not introduced or even mentioned.

\section{A synthetic approach}

In the schools where the synthetic definition of a parabola is also
introduced, the most common definition is that it is \textit{the locus of points in the plane that
are equidistant from both the directrix line $\ell$ and the focus point $F$}.

\subsection{An automated answer}

\begin{figure}
\begin{center}
\includegraphics[scale=0.4]{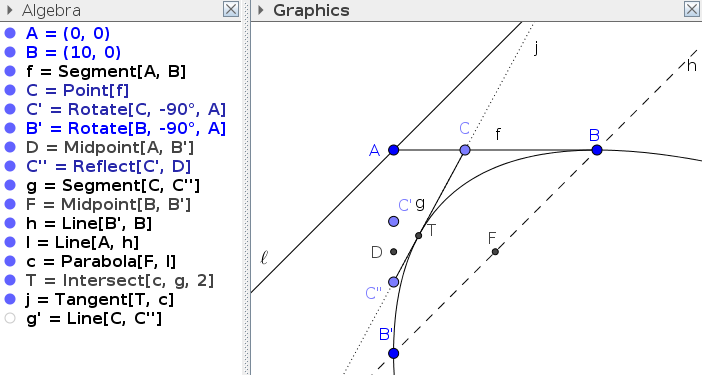}
\caption{A synthetic way to construct the investigated parabola in GeoGebra}
\label{fig:parabola4}
\end{center}
\end{figure}

Without any further considerations it is possible to check (actually, \textit{prove})
that the investigated curve is a parabola also in the symbolic sense. To achieve
this result one can invoke GeoGebra's Relation Tool \cite{RelTool-ADG2014} after constructing the parabola
synthetically as seen in Fig.~\ref{fig:parabola4}.

Here the focus point $F$ is the midpoint of $BB'$, and the directrix line $\ell$
is a parallel line to $BB'$ through $A$. This piece of information should be
probably kept secret by the teacher---the learners could find them on their
own. Now to check if the string $g$ is indeed a tangent of the parabola the
tangent point $T$ has been created as an intersection of the string and
the parabola. Also a tangent line $j$ has been drawn, and finally line $g'$
which is the extension of segment $g$ to a full line. At this point it is
possible to compare $g'$ and $j$ by using the Relation Tool.

\begin{figure}
\begin{center}
\includegraphics[scale=0.4]{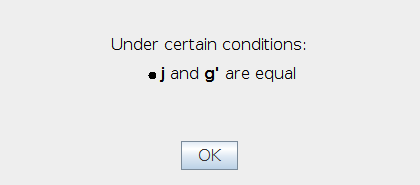}
\caption{Symbolic check of the synthetic statement in GeoGebra 5.0.354.0.
Here ``under certain conditions'' mean that the statement holds in general,
but there may be some extra prerequisites like avoiding degenerate positions
of the free points, which cannot be further described by GeoGebra \cite{gg-art-doc-gh}.}
\label{fig:relation}
\end{center}
\end{figure}

The Relation Tool first compares the two objects numerically and reports
that they are equal. By clicking the ``More$\ldots$'' button the user
obtains the symbolic result of the synthetic statement (Fig.~\ref{fig:relation}).

We recall that despite the construction was performed synthetically, the symbolic computations
were done after translating the construction to an algebraic setup.
Thus GeoGebra's internal proof is again based on algebraic equations and
still hidden for the user. But in this case we indeed have a general proof for
each possible construction setup, not for only one particular case as
for the \texttt{Envelope} command.

\subsection{A classical proof}

Finally we give a classical proof to answer the original question. Here every detail uses
only synthetic considerations.

The first part of the proof is a well known remark on the bisection
property of the tangent. That is, by reflecting the focus point about any
tangent of the parabola the mirror image is a point of the directrix line.
(See e.g.~\cite{LociEnvelopes-arXiv}, Sect.~3.1 for a short proof.) Clearly,
it is sufficient to show that the strings have this kind of bisection property:
this will result in confirming the statement.

\begin{figure}
\begin{center}
\includegraphics[scale=0.5]{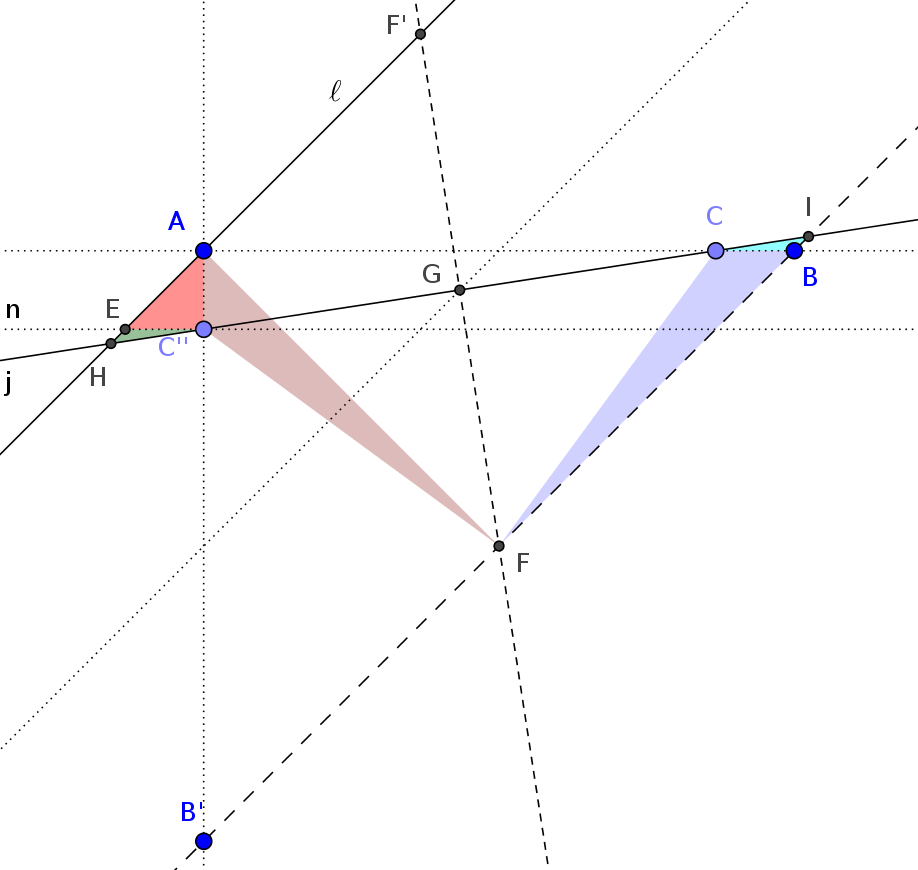}
\caption{A classical proof}
\label{fig:parabola5}
\end{center}
\end{figure}

In Fig.~\ref{fig:parabola5} the tangent to the parabola is denoted by $j$.
Let $F'$ be the mirror image of $F$ about $j$. We will prove that $F'\in\ell$.
Let $G$ denote the intersection of $j$ and $FF'$. Clearly $\angle CGF$ and
$\angle C''GF$ are right because of the reflection.

By construction $\triangle FBC$ and $\triangle FAC''$ are congruent.
Thus $FC=FC''$. Moreover, $\triangle CFC''$ is isosceles and $FG$ is
its bisector at $F$, in addition, $CG=C''G$.

Let $n$ denote a parallel line to $AB$ through $C''$. Let $E$ be the intersection of $\ell$ and $n$.
Also let $H$ be the intersection of $j$ and $\ell$, and let $I$ be the
intersection of $j$ and the line $BB'$. Since $\ell$ and $BB'$ are parallel,
moreover $AB$ and $n$ are also parallel, and $EC''=CB$ (because $E$ is actually
the rotation of $A$ around $C''$ by 90 degrees), we conclude that $\triangle C''EH$
and $\triangle CBI$ are congruent. This means that $IC=C''H$.

That is, using also $CG=C''G$, $G$ must be the midpoint of $HI$, thus $G$ lies on the mid-parallel of $\ell$ and $BB'$.
As a consequence, reflecting $F$ about $G$ the resulted point $F'$ is surely
a point of line $\ell$.

\section{Conclusion}

An analysis of the string art envelope was presented at different levels of
mathematical knowledge, by refuting a false conjecture, finding a true
statement and then proving it with various means.

Discussion of a non-trivial question by using different means can give
a better understanding of the problem. What is more, reasoning by visual ``evidence'' can
be misleading, and only rigorous (or rigorous \textit{but} computer based) proofs
can be satisfactory.

\begin{figure}
\begin{center}
\includegraphics[scale=0.3]{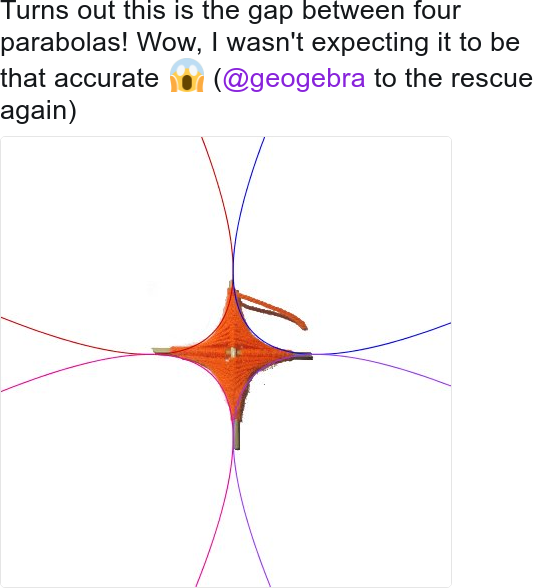}
\caption{Another way to experiment with the string art parabola \cite{tweet:geneeper}}
\label{fig:tweet}
\end{center}
\end{figure}

It should be noted that the parabola property of the string art envelope
is well known in the literature on Bezier curves, but usually not
discussed in maths teacher trainings. The \textit{de Casteljeau
algorithm} for a Bezier curve of degree 2 is itself a proof that the
curve is a parabola. (See \cite{wiki:parabola,Markus-stringart} for more
details.) Also among maths professionals this property seems rarely
known. A recent example of a tweet of excitement is from February 2017
(Fig. \ref{fig:tweet}).

On the other hand, our approach highlighted the classroom introduction
of the string art parabola, and suggested some very recent methods by
utilizing computers in the middle and high school to improve the teacher's work
and the learners' skills.

Lastly, we remark that the definition of the string art envelope looks
similar to the envelope of other family of lines. For example, the
envelope of the sliding ladder results in a different curve, the
\textit{astroid} \cite{wiki:envelope,wiki:astroid,sliding-ladder-ggm}, a
real algebraic curve of degree 6. While ``physically'' that is easier to
construct (one just needs a ladder-like object, e.~g.~a pen), the geometric
analysis of that is more complicated and usually involves partial
derivatives. (See also \cite{wiki:envelope} on a proof for identifying the string art
parabola by using partial derivatives.)

\section*{Acknowledgments}
The author thanks Tom\'as Recio and No\"el Lambert for comments that greatly improved the manuscript.

The introductory figure was drawn by Benedek Kov\'acs (12), first grade middle school student.

\bibliography{kovzol,external}

\end{document}